\documentclass[12pt]{amsart}
\usepackage{amssymb, mathrsfs}

\newcommand{\ed}{

\end{document}
}

\long\def\forget#1\forgotten{}
\newcommand{\Pa}[9]{\bibitem{#1} {#2}, \emph{#3}, {#4} \textbf{#5} ({#6}), {#7}--{#8}.}
\newcommand{\Bc}[9]{\bibitem{#1} {#2}, \emph{#3}, in: \textbf{#4} (#5), #6 #7, #8--#9.}

\newcommand{\op}[1]{\operatorname{#1}}
\newcommand{\Par}{\op{Par}}

\newcommand{\AP}{\mathrm{A\!P}}

\newcommand{\Impl}{\Rightarrow}
\newcommand{\bi}{\begin{itemize}}
\newcommand{\itm}{\item}
\newcommand{\ei}{\end{itemize}}
\newcommand{\be}{\begin{enumerate}}
\newcommand{\ee}{\end{enumerate}}

\newcommand{\nin}{\notin}
\newcommand{\cO}{\mathcal{O}}
\newcommand{\Union}{\bigcup}
\newcommand{\sbst}{\subseteq}
\newcommand{\spst}{\supseteq}
\newcommand{\N}{\mathbb{N}}

\newcommand{\bbB}{\mathbb{B}}

\newcommand{\seq}[1]{\{#1\}_{n\in\N}}
\newcommand{\sseq}[1]{\{#1 : n\in\N\}}
\newcommand{\Icov}[1]{{#1}_{\cI}}

\newcommand{\sm}{\setminus}
\newcommand{\FinSeqs}[1]{{#1}^*}

\newcommand{\roth}{{[\N]^{\!\/\infty}}}
\newcommand{\setroth}[1]{{[#1]^{\!\/\infty}}}

\newcommand{\sone}{\mathsf{S}_{1}}
\newcommand{\sfin}{\mathsf{S}_\mathrm{fin}}

\newcommand{\cA}{\mathcal{A}}
\newcommand{\cI}{\mathcal{I}}
\newcommand{\OI}{\cO_\cI}

\newcommand{\cU}{\mathcal{U}}
\newcommand{\cV}{\mathcal{V}}

\newcommand{\cF}{\mathcal{F}}

\newcommand{\cS}{\mathcal{S}}

\newtheorem{thm}{Theorem}[section]
\newcommand{\bthm}{\begin{thm}} \newcommand{\ethm}{\end{thm}}
\newtheorem{prop}[thm]{Proposition}
\newcommand{\bprp}{\begin{prop}} \newcommand{\eprp}{\end{prop}}
\newtheorem{fact}[thm]{Fact}
\newcommand{\bfct}{\begin{fact}} \newcommand{\efct}{\end{fact}}
\newtheorem{prob}[thm]{Problem}
\newcommand{\bprb}{\begin{prob}} \newcommand{\eprb}{\end{prob}}
\newtheorem{lem}[thm]{Lemma}
\newcommand{\blem}{\begin{lem}} \newcommand{\elem}{\end{lem}}
\newtheorem{claim}[thm]{Claim}
\newcommand{\bclm}{\begin{claim}} \newcommand{\eclm}{\end{claim}}
\newtheorem{cor}[thm]{Corollary}
\newcommand{\bcor}{\begin{cor}} \newcommand{\ecor}{\end{cor}}
\newtheorem{conj}[thm]{Conjecture}
\newcommand{\bcnj}{\begin{conj}} \newcommand{\ecnj}{\end{conj}}
\theoremstyle{definition}
\newtheorem{defn}[thm]{Definition}
\newcommand{\bdfn}{\begin{defn}} \newcommand{\edfn}{\end{defn}}
\theoremstyle{remark}
\newtheorem{rem}[thm]{Remark}
\newcommand{\brem}{\begin{rem}} \newcommand{\erem}{\end{rem}}
\newtheorem{cnv}[thm]{Convention}
\newcommand{\bcnv}{\begin{cnv}} \newcommand{\ecnv}{\end{cnv}}
\newtheorem{exam}[thm]{Example}
\newcommand{\bexm}{\begin{exam}} \newcommand{\eexm}{\end{exam}}
\newcommand{\bpf}{\begin{proof}} \newcommand{\epf}{\end{proof}}

\title[Superfilters and Ramsey theory]{Superfilters, Ramsey theory, and van der Waerden's Theorem}

\author{Nadav Samet}
\address[Nadav Samet]{Department of Mathematics,
Weizmann Institute of Science, Rehovot 76100, Israel
}
\curraddr{Google Ireland Ltd., Gordon House, Barrow Street, Dublin 4, Ireland}
\email{thesamet@gmail.com}

\author{Boaz Tsaban}
\address[Boaz Tsaban]{Department of Mathematics, Bar-Ilan University, Ramat-Gan 52900, Israel;
and
Department of Mathematics, Weizmann Institute of Science, Rehovot 76100, Israel
}
\email{tsaban@math.biu.ac.il}

\begin{document}
\begin{abstract}
Superfilters are generalizations of ultrafilters, and capture the
underlying concept in Ramsey theoretic theorems such as van der
Waerden's Theorem. We establish several properties of
superfilters, which generalize both Ramsey's Theorem and its
variants for ultrafilters on the natural numbers. We use them to
confirm a conjecture of Ko\v{c}inac and Di Maio, which is a
generalization of a Ramsey theoretic result of Scheepers,
concerning selections from open covers. Following Bergelson and
Hindman's 1989 Theorem, we present a new simultaneous
generalization of the theorems of Ramsey, van der Waerden, Schur,
Folkman-Rado-Sanders, Rado, and others, where the colored sets
can be much smaller than the full set of natural numbers.
\end{abstract}

\maketitle

\section{A unified Ramsey Theorem}

It is a simple observation that when each element of an infinite set is
colored by one of finitely many colors,
the set must contain an infinite monochromatic subset.
When replacing \emph{infinite} by \emph{containing arithmetic progressions of arbitrary length},
we obtain van der Waerden's Theorem \cite{vdW}.
Some of the best references for many beautiful theorems of this kind, together with applications,
are the classical \cite{GRS}, the monumental \cite{HS}, the elegant
Protasov \cite{ProtasovRT}, and the more recent \cite{LR}.
These results lead naturally to the concept of superfilter.

\bdfn\label{suf}
For a set $S$, $[S]^n=\{F\sbst S : |F|=n\}$, and $\setroth{S}$
is the family of infinite subsets of $S$.

A nonempty family $\cS\sbst\roth$ is a \emph{superfilter} if for all $A,B\sbst\N$:
\be
\itm If $A\in\cS$ and $B\spst A$, then $B\in\cS$.
\itm If $A\cup B\in\cS$, then $A\in\cS$ or $B\in\cS$.
\ee
\edfn

Superfilters were identified at least as early as in Berge's 1959
monograph \cite{Berge} (under the name \emph{grille}).\footnote{(Added after publication)
Frederic Mynard points out that the notion of superfilter (under the name of grille) goes back at least to:
Gustave Choquet, \emph{Sur les notions de filtre et de grille}, C. R. Acad.\ Sci.\ Paris \textbf{224} (1947), 171--173.}
They were
also considered under the name \emph{coideal} (e.g.,
\cite{Farah}). Superfilters are large types of Banakh and
Zdomskyy's \emph{semifilters} and \emph{unsplit} semifilters \cite{CSF}.

Recall that a \emph{nonprincipal ultrafilter} is a family as in Definition \ref{suf} which is also closed
under finite intersections.\footnote{
Definition \ref{suf} does not change if we assume that $A,B$ are disjoint in (2).
But if, in addition, we replace there \emph{or} by \emph{exclusive or}, we obtain a characterization
of ultrafilter. That is, the assumption about intersections need not be stated
explicitly.}
For brevity, by \emph{ultrafilter} we always mean a nonprincipal one.

\bexm
\mbox{}
\be
\itm Every ultrafilter is a superfilter.
\itm Every union of a family of ultrafilters is a superfilter.
\itm $\roth$ is a superfilter which is not an ultrafilter.
\ee
\eexm
In fact, one can show that every superfilter is a union of a family of
ultrafilters, but we will not use this here.

\bdfn
$\AP$ is the family of all subsets of $\N$ containing arbitrarily long arithmetic progressions.
\edfn

Clearly, $\AP$ is not an ultrafilter.
The finitary version of van der Waerden's Theorem implies the following.

\bthm[van der Waerden]\label{vdW}
$\AP$ is a superfilter.
\ethm

\forget
\brem
For each superfilter $\cS$:
\be
\itm For each $A\in\cS$ and each finite decomposition $A=A_1\cup A_2\cup\dots\cup A_k$,
there is $i\le k$ such that $A_i\in\cS$.
\itm If $A_1,\dots,A_k\nin\cS$, then $A_1\cup\dots\cup A_k\nin\cS$.
\itm For each $A\in\cS$ and each $B\nin\cS$, $A\sm B\in\cS$.
\ee
\erem
\bpf
(1) Induction on $k$.

(2) This is the contrapositive of (1).

(3) $(A\sm B)\cup(A\cap B)=A\in\cS$.
\epf
\forgotten

\bdfn\label{Ram}
$\cS\to(\cS)^n_k$ is the statement: For each $A\in\cS$ and each coloring $c:[A]^n\to\{1,2,\dots,k\}$,
there is $M\sbst A$ such that $M\in\cS$ and $c$ is constant on $[M]^n$.
The set $M$ is called \emph{monochromatic} for the coloring $c$.
\edfn

Thus for upwards-closed $\cS\sbst\roth$, the following are equivalent:
\be
\itm $\cS$ is a superfilter.
\itm $\cS\to(\cS)^1_2$.
\itm $\cS\to(\cS)^1_k$ for all $k$.
\ee
The assertion $\cS\to(\cS)^n_k$ becomes stronger when $n$ or $k$ is increased.

\bdfn
A superfilter $\cS$ is:
\be
\itm \emph{Ramsey} if $\cS\to(\cS)^n_k$ holds for all $n$ and $k$.
\itm \emph{Strongly Ramsey} if for all pairwise disjoint $A_1,A_2,\dots$ with
$\Union_{n\ge m} A_n\in\cS$ for all $m$, there is $A\sbst\Union_n A_n$ such that $A\in\cS$ and
$|A\cap A_n|\le 1$ for all $n$.
\itm \emph{Weakly Ramsey} if for all pairwise disjoint $A_1,A_2,\dots\nin\cS$ with
$\Union_n A_n\allowbreak\in\cS$, there is $A\sbst\Union_n A_n$ such that $A\in\cS$ and
$|A\cap A_n|\le 1$ for all $n$.
\ee
\edfn

Clearly, strongly Ramsey superfilters are weakly Ramsey. We will
soon show that Ramsey is sandwiched between strongly Ramsey and
weakly Ramsey. Before doing so, we give examples showing that
converse implications cannot be proved.

\bexm\label{RnotsR}
Fix a partition $\N=\Union_nI_n$ with each $I_n$ infinite.
Let $\cS$ be the upwards closure of $\Union_n\setroth{I_n}$.
It is easy to see that $\cS$ is a superfilter.

$\cS$ is Ramsey:
Let $A\in\cS$, and $c:[A]^n\to\{1,\dots,k\}$ be a coloring of $A$. Pick $m$ such that
$A\cap I_m$ is infinite, and use Ramsey's Theorem \ref{Ramsey} for the coloring
$c:[A\cap I_m]^2\to\{1,\dots,k\}$ to obtain an infinite $M\sbst A\cap I_m$ which is
monochromatic for $c$.

$\cS$ is not strongly Ramsey: For each $m$, $\Union_{n\ge m} I_n\in\cS$, but if
$|A\cap I_n|\le 1$ for all $n$, then $A\nin\cS$.
\eexm

\bexm\label{wRnotR}
Following is an example of a weakly Ramsey superfilter which is not Ramsey.
Essentially the same example was, independently, found by Filip\'ow, Mro\.zek, Rec\l{}aw,
and Szuca \cite{FMRS}.

Let $\FinSeqs\N$ be the set of all finite sequences of natural numbers.
For $\sigma,\rho\in\FinSeqs\N$, write $\sigma\spst\rho$ if the sequence
$\rho$ is a prefix of $\sigma$.
As $\FinSeqs\N$ is countable, we may use it instead of $\N$ to define our superfilter.
Say that a set $D\sbst\FinSeqs\N$ is \emph{somewhere dense} if there is $\rho\in\FinSeqs\N$
such that for each $\sigma\in\FinSeqs\N$ with $\sigma\spst\rho$, there is
$\eta\spst\sigma$ such that $\eta\in D$.
Let $\cS$ be the family of all somewhere dense subsets of $\FinSeqs\N$.

It is not difficult to see that $\cS$ is a superfilter, and that it is weakly Ramsey.
To see that it is not Ramsey, define a coloring $c:[\FinSeqs\N]^2\to\{1,2\}$ by
$c(\sigma,\eta)=1$ if one of $\sigma,\eta$ is a prefix of the other, and $2$ otherwise.
If $M\sbst\FinSeqs\N$ is monochromatic of color $1$, then $M$ is a branch
in $\FinSeqs\N$, and thus $M\nin\cS$. On the other hand, if $M$ is
somewhere dense, then it must contain at least two elements, one of which a prefix of
the other. Thus, $M$ is not monochromatic of color $2$, either.
\eexm

Examples \ref{RnotsR} and \ref{wRnotR} show that some hypothesis is
required to make the Ramseyan notions coincide. We suggest a rather
mild one.

\bdfn
A superfilter $\cS$ is \emph{shrinkable} if,
for all pairwise disjoint $A_1,A_2,\dots$ with $\Union_{n\ge m} A_n\in\cS$ for all $m$,
there are $B_n\sbst A_n$ such that $B_n\nin\cS$ and $\Union_nB_n\in\cS$.
\edfn

\brem[Thuemmel]
A superfilter $\cS$ is shrinkable if, and only if,
for each sequence $S_1\spst S_2\spst\dots$ of element of $\cS$,
there is $S\in\cS$ such that for each $n$, $S\sm S_n\nin\cS$.
To see this, identify $S_m$ with $\Union_{n\ge m} A_n$ for each $m\in\N$,
and $S$ with $\Union_nB_n$.
\erem

All ultrafilters are shrinkable, for a trivial reason:
If a disjoint union $\Union_nA_n$ is in the ultrafilter, and some
$A_m$ is in the ultrafilter, then $\Union_{n>m}A_n$ is not in the ultrafilter.

The superfilters in Examples \ref{RnotsR} and \ref{wRnotR} are not shrinkable.
For shrinkable superfilters, we have a complete characterization of being Ramsey.

\bthm\label{main}
For superfilters $\cS$, the following are equivalent:
\be
\itm $\cS$ is strongly Ramsey.
\itm $\cS$ is Ramsey and shrinkable.
\itm $\cS\to(\cS)^2_2$, and $\cS$ is shrinkable.
\itm $\cS$ is weakly Ramsey and shrinkable.
\ee
\ethm
\bpf
$(1\Impl 2)$ As singletons do not belong to superfilters, strongly Ramsey implies shrinkable.
It therefore suffices to prove the following.

\blem\label{sRisR}
Every strongly Ramsey superfilter is Ramsey.
\elem
\bpf
Let $\cS$ be a strongly Ramsey superfilter,
$A\in\cS$, and $c:[A]^d\to\{1,\dots,k\}$.
The proof is by induction on $d$, with $d=1$
following from $\cS$ being a superfilter.

Induction step: We repeatedly apply the following fact.
For each $A\in\cS$ and each $n\in A$,
there is $M\sbst A\sm\{n\}$ such that $M\in\cS$, and
a color $i\in\{1,\dots,k\}$, such that for each $F\in[M]^{d-1}$,
$c(\{n\}\cup F)=i$. Indeed, we can define a coloring
$c_n:[A\sm\{n\}]^{d-1}\to\{1,\dots,k\}$ by $c_n(F)=c(\{n\}\cup F)$
and use the induction hypothesis.

Enumerate $A=\sseq{a_n}$.
Choose $A_{a_1}\sbst A\sm\{a_1\}$ and a color $i_{a_1}$ such that
$A_{a_1}\in\cS$ and for each $F\in[A_{a_1}]^{d-1}$, $c(\{a_1\}\cup F)=i_{a_1}$.
In a similar manner, choose inductively for each $n>1$
$A_{a_n}\sbst A_{a_{n-1}}\sm\{a_n\}$ and a color $i_{a_n}$ such that
$A_{a_n}\in\cS$ and for each $F\in[A_{a_n}]^{d-1}$, $c(\{a_n\}\cup F)=i_{a_n}$.

As $a_n\nin A_{a_n}$ for all $n$, $\bigcap_n A_{a_n}=\emptyset$.
Let $B_0=A\sm A_{a_1}$ and for each $n>0$, let $B_n=A_{a_n}\sm A_{a_{n+1}}$.
The sets $B_n$ are pairwise disjoint, $\Union_nB_n=A$, and
$\Union_{n\ge m}B_n=A_{a_m}\in\cS$ for all $m$.
As $\cS$ is strongly Ramsey, there is $B\sbst A$ such that $B\in\cS$ and
$|B\cap B_n|\le 1$ for all $n$.
Fix a color $i$ such that $C=\{n\in B : i_n=i\}\in\cS$.

Let $c_1=\min C$.
Inductively, for each $n>1$ choose $c_n\in C$ such that $c_n>c_{n-1}$
and $C\sm[1,c_n)\sbst A_{c_{n-1}}$.\footnote{E.g., let $k=|C\sm A_{c_{n-1}}|+1$ and let $c_n$ be the
$k$-th element of $C$.}
For each $n$, $C\cap[c_n,c_{n+1})$ is finite and thus not a member of $\cS$.
As $\Union_n(C\cap[c_n,c_{n+1}))=C\in\cS$ and $\cS$ is weakly Ramsey,
there is $D\in\cS$ such that $D\sbst C$ and $|D\cap [c_n,c_{n+1})|\le 1$ for all $n$.
As
$$D=\left(D\cap \Union_{n\in\N} [c_{2n},c_{2n+1})\right)\cup\left(D\cap \Union_{n\in\N}[c_{2n-1},c_{2n})\right),$$
there is $l\in\{0,1\}$ such that $M=D\cap \Union_n[c_{2n-l},c_{2n+1-l})\in\cS$.
Let $m_1<m_2<\dots<m_d$ be members of $M$. Let $n$ be minimal such that $m_1<c_n$. Then
$$m_2,\dots,m_d\in C\sm[1,c_{n+1})\sbst A_{c_n}\sbst A_{m_1},$$
and thus $c(\{m_1,\dots,m_d\})=c(\{m_1\}\cup\{m_2,\dots,m_d\}) = i_{m_1}=i$.
\epf

$(2\Impl 3)$ Trivial.

$(3\Impl 4)$ In fact, the following holds.

\blem\label{RiswR}
If $\cS\to(\cS)^2_2$, then $\cS$ is weakly Ramsey.
\elem
\bpf
Let $A_1,A_2,\dots$ be as in the definition of weakly Ramsey.
Let $D=\Union_n A_n$, and define a coloring $c:[D]^2\to\{1,2\}$ by
$$c(m,k)=\begin{cases}
1 & (\exists n)\ m,k\in A_n\\
2 & \mbox{otherwise}
\end{cases}$$
As $\cS$ is Ramsey, there is a monochromatic $A\sbst D$ with $A\in\cS$.
If all elements of $[A]^2$ have color $1$, then $A\sbst A_n$ for some $n$,
and thus $A_n\in\cS$, a contradiction.
Thus, all elements of $[A]^2$ have color $2$, which means that $|A\cap A_n|\le 1$ for all $n$.
\epf

$(4\Impl 1)$ Let $A_1,A_2,\dots$ be as in the definition of strongly Ramsey.
As $\cS$ is shrinkable, there are $B_n\sbst A_n$ such that $B_n\nin\cS$ and $B=\Union_nB_n\in\cS$.
As $\cS$ is weakly Ramsey, there is a subset $A$ of $B$ such that $A\in\cS$ and $|A\cap B_n|\le 1$ for all $n$.
As $B_n\sbst A_n$ for all $n$ and the sets $A_n$ are pairwise disjoint, $|A\cap A_n|\le 1$ for all $n$.

This completes the proof of Theorem \ref{main}.
\epf

\bcor[Ramsey \cite{Ramsey}]\label{Ramsey}
$\roth\to(\roth)^n_k$ for all $n$ and $k$.
\ecor
\bpf
Clearly, $\roth$ is strongly Ramsey.
\epf

\bcor[Booth-Kunen \cite{Booth}]\label{ufR}
An ultrafilter is weakly Ramsey if, and only if, it is Ramsey.
\ecor
\bpf
Ultrafilters are shrinkable.
\epf

The following definition and subsequent result will be useful later.

\bdfn[Scheepers \cite{coc1}]
$\sone(\cS,\cS)$ is the statement: Whenever $S_1,S_2,\dots\in\cS$, there
are $s_n\in S_n$, $n\in\N$, such that $\sseq{s_n}\in\cS$.
\edfn

\bthm\label{shrink}
For superfilters $\cS$:
\be
\itm If $\cS$ is strongly Ramsey, then $\sone(\cS,\cS)$ holds.
\itm $\sone(\cS,\cS)$ implies that $\cS$ is shrinkable.
\ee
\ethm
\bpf
(1) We first observe that, in the definition of strongly Ramsey,
there is no need for the sets $A_n$ to be pairwise disjoint.
\blem\label{nodisj}
If a superfilter $\cS$ is strongly Ramsey, then
for all nonempty $A_1,A_2,\dots$ with $\Union_{n\ge m} A_n\in\cS$ for all $m$,
there are $a_n\in A_n$, $n\in\N$,
such that $A=\sseq{a_n}\in\cS$.
\elem
\bpf
Assume that $\Union_{n\ge m} A_n\in\cS$ for all $m$. Let
$$L=\bigcap_{m\in\N}\bigcup_{n\ge m}A_n.$$
If $L\in\cS$, enumerate $L=\sseq{l_n}$.
Pick $m_1$ such that $a_{m_1}:=l_1\in A_{m_1}$.
For each $n>1$, there is $m_n>m_{n-1}$ such that $a_{m_n}:=l_n\in A_{m_n}$.
For $m\nin\seq{m_n}$, pick any $a_m\in A_m$.
Then we obtain a sequence as required.

Thus, assume that $L\nin\cS$.
Taking $B_n=A_n\sm L$ for all $n$, we have that
$$\Union_{n\ge m} B_n=(\Union_{n\ge m}A_n)\sm L\in\cS$$
for all $m$.
Now, $\bigcap_m\bigcup_{n\ge m}B_n=\emptyset$, that is, each $k\in\Union_nB_n$
belongs to only finitely many $B_n$.
For each $n$, let
$$C_n = B_n\sm\Union_{m>n}B_m.$$
The sets $C_n$ are pairwise disjoint, and for each $m$,
$\Union_{n\ge m}C_n=\Union_{n\ge m}B_n\allowbreak\in\cS$.
As $\cS$ is strongly Ramsey, we obtain $A\sbst\Union_nC_n$ such
that $A\in\cS$ and $|A\cap C_n|\le 1$ for all $n$.
For each $n$, let $a_n\in A\cap C_n$ if $|A\cap C_n|=1$, and
an arbitrary element of $A_n$ otherwise.
Then the sequence $\seq{a_n}$ is as required.
\epf
Thus, assume that $A_1,A_2,\dots\in\cS$. Clearly, they are all nonempty,
and $\Union_{n\ge m}A_n\in\cS$ for all $m$.
By Lemma \ref{nodisj}, there are $a_n\in A_n$, $n\in\N$, such that $\sseq{a_n}\in\cS$.

(2) Apply $\sone(\cS,\cS)$ to the sequence $\Union_{n\ge m} A_n$, $m\in\N$,
and recall that finite sets do not belong to superfilters.
\epf

As Ramsey does not imply strongly Ramsey (Example \ref{RnotsR}), but does
for shrinkable superfilters (Theorem \ref{main}(4)), we have that the
converse of Theorem \ref{shrink}(2) is false.
Unfortunately, we do not have a concrete example for the following.

\bcnj
There is a superfilter $\cS$ such that $\sone(\cS,\cS)$ holds, but
$\cS$ is not strongly (equivalently, by Theorem \ref{shrink}(2), weakly) Ramsey.
\ecnj

\section{An application to topological selection principles}

Our initial motivation for studying superfilters came from an attempt to
provide a (mainly) combinatorial proof of a major Ramsey-theoretic result of Scheepers,
concerning selections from open covers. The general theory has connections and applications
far beyond Ramsey theory, and the interested reader is referred to the survey papers
\cite{LecceSurvey, KocSurv, ict}. The Ramsey-theoretic aspect of this theory is surveyed
in \cite{KocRamsey}.
Here, we present only the concepts which are necessary for the present paper.

Fix a topological space $X$.
A family $\cU$ of subsets of $X$ is a \emph{cover} of $X$ if $X\nin\cU$ but $X=\Union\cU$.
A cover $\cU$ of $X$ is an \emph{$\omega$-cover} if for each finite $F\sbst X$,
there is $U\in\cU$ such that $F\subseteq U$.
Let $\Omega=\Omega(X)$ denote the family of all open $\omega$-covers of $X$.
According to Definition \ref{Ram}, the statement $\Omega\to(\Omega)^2_2$ makes sense,
and it is natural to ask what is required from $X$ for this
statement to be true.
Say that $X$ is \emph{$\Omega$-Lindel\"of} if each element of $\Omega$ contains a countable element of $\Omega$.
The following result is essentially proved in \cite{coc1}, using an auxiliary result from \cite{coc2}.
In the general form stated here, it is proved in \cite{coc7}.

\bthm[Scheepers \cite{coc1, coc2, coc7}]\label{sch}
For $\Omega$-Lindel\"of spaces, the following are equivalent:
\be
\itm $\sone(\Omega,\Omega)$.
\itm $\Omega\to(\Omega)^2_2$.
\itm $\Omega\to(\Omega)^n_k$ for all $n,k$.
\ee
\ethm

We proceed in a general manner that will prove, in addition to
Scheepers's Theorem, a conjecture of Di Maio, Ko\v{c}inac, and Meccariello from \cite{AppKI},
and a subsequent one of Di Maio and Ko\v{c}inac from \cite{DMKBdd}.

Let $C(X)$ denote the space of continuous real-valued functions of $X$.
$\omega$-covers arise when considering the closure operator in $C(X)$,
with the topology of pointwise convergence \cite{GN}. When considering the
compact-open topology, \emph{$k$-covers} arise, which are covers such that each
compact set is contained in a member of the cover (e.g., \cite{AppKI}
and references therein). In \cite{AppKI} it is conjectured that Scheepers's
Theorem also holds when $\omega$-covers are replaced by $k$-covers.

A natural generalization of these topologies
on $C(X)$ gives rise to the following notion.
An \emph{abstract boundedness} is a family $\bbB$ of nonempty closed subsets of $X$ which
is closed under taking finite unions and closed subsets, and contains all singletons
\cite{DMKBdd}.
A cover $\cU$ is a \emph{$\bbB$-cover} if each $B\in\bbB$ is
contained in some member of $\cU$.
In \cite{DMKBdd} it is conjectured that Scheepers's
Theorem holds in general, when $\omega$-covers are replaced by $\bbB$-covers for
any abstract boundedness notion $\bbB$.

Closing an abstract boundedness notion $\bbB$ downwards will not
change the notion of $\bbB$-covers. Thus, for simplicity we use a more familiar
notion.
A nonempty family $\cI$ of subsets of $X$ is an \emph{ideal} on $X$ if
$X\nin\cI$, $\{x\}\in\cI$ for all $x\in X$,
and for all $A,B\in\cI$, $A\cup B\in\cI$.

\bdfn\label{afterideal}
Fix an ideal $\cI$ on $X$.
$\cU$ is an $\cI$-cover of $X$ if $X\nin\cU$, and for each $B\in\cI$ there is
$U\in\cU$ such that $B\sbst U$. $\OI$ is the family of all open $\cI$-covers of $X$.
\edfn

\blem\label{folk}
\mbox{}
\be
\itm If $\cU_1\cup\cU_2\in\OI$, then $\cU_1\in\OI$ or $\cU_2\in\OI$.
\itm Each $\cU\in\OI$ is infinite.
\ee
\elem
\bpf
(1) Assume that $B_1,B_2\in\cI$ witness that $\cU_1,\cU_2\nin\OI$, respectively.
Then no element of $\cU_1\cup\cU_2$ contains $B_1\cup B_2$.

(2) $\OI\sbst\Omega$.
\epf

Let $\cU\in\OI$. If $\cU$ is countable, we may use it as an index set instead of $\N$,
and consider superfilters on $\cU$.

\bdfn
$\Icov\cU=\{\cV\sbst\cU : \cV\in\OI\}=P(\cU)\cap\OI$.
\edfn

Lemma \ref{folk} implies the following.

\bcor\label{wsuf}
For each countable $\cU\in\OI$, $\Icov\cU$ is a superfilter.\qed
\ecor

$\Icov\cU$ cannot be assumed to be an ultrafilter when proving Scheepers's Theorem \ref{sch}:
If $\sone(\Omega,\Omega)$ holds, then each $\cU\in\Omega$ can be split into two disjoint elements of
$\Omega$ \cite{coc1}.

We are now ready to prove the general statement.
Say that $X$ is \emph{$\OI$-Lindel\"of} if each element of $\OI$ contains a countable element of $\OI$.

\bthm\label{sch++}
Let $\cI$ be an ideal on $X$.
For $\OI$-Lindel\"of spaces, the following are equivalent:
\be
\itm $\sone(\OI,\OI)$.
\itm For all disjoint $\cU_1,\cU_2,\dots\nin\OI$ with $\Union_n\cU_n\in\OI$,
there is $\cV\sbst\Union_n\cU_n$ such that $\cV\in\OI$ and $|\cV\cap\cU_n|\le 1$ for all $n$.
\itm $\OI\to(\OI)^2_2$.
\itm $\OI\to(\OI)^n_k$ for all $n,k$.
\ee
\ethm
\bpf
Using $\OI$-Lindel\"ofness, we may restrict attention to countable $\cI$-covers
in all of our arguments. More precisely, we prove the stronger assertion, where
$\OI$ is replaced with the family of \emph{countable} open $\cI$-covers,
and no assumption is posed on the space $X$.

$(4\Impl 3)$ Trivial.

$(3\Impl 2)$ Let $\cU_1,\cU_2,\dots$ be as in $(2)$. Set $\cU=\Union_n\cU_n$.
Then $\cU\in\OI$, and by Corollary \ref{wsuf}, $\Icov\cU$ is a superfilter.
By $(3)$, we have in particular $\Icov\cU\to(\Icov\cU)^2_2$.
By Theorem \ref{RiswR}, $\Icov\cU$ is weakly Ramsey.
As $\cU_1,\cU_2,\dots\nin\Icov\cU$ and $\Union_n\cU_n=\cU\in\Icov\cU$,
there is $\cV\in\Icov\cU\sbst\OI$ as required.

$(2\Impl 1)$ Assume that $\cU_1,\cU_2,\dots\in\OI$.
Fix $\sseq{U_n}\in\OI$.
For each $n$, let
$$\cV_n=\{U_n\cap U : U\in\cU_n\}.$$
Then
$$\cU=\Union_{n\in\N}\cV_n\in\OI.$$
By Corollary \ref{wsuf}, $\Icov\cU$ is a superfilter.
By $(2)$, $\Icov\cU$ is weakly Ramsey.
Now, $\Union_n\cV_n=\cU\in\Icov\cU$, and for each $n$, $\cV_n\nin\Icov\cU$.
By thinning out the sets $\cV_n$ if necessary, we may assume that they are disjoint.
Thus, there is $\cV\sbst\cU$ such that $\cV\in\Icov\cU$
and $|\cV\cap\cV_n|\le 1$ for all $n$.

For each $n$, if $|\cV\cap\cV_n| = 1$, take the $U\in\cU_n$
such that $U_n\cap U\in\cV$, and otherwise take an arbitrary $U\in\cU_n$.
We obtain a $\cI$-cover of $X$ with one element from each $\cU_n$.

$(1\Impl 4)$ Let $\cU\in\OI$.
Let $\cV$ be the closure of $\cU$ under finite intersections.
$\cV$ is countable, and $\cU\in\Icov\cV\sbst\OI$.

Consider the superfilter $\Icov\cV$.
By $\sone(\OI,\OI)$, we have $\sone(\Icov\cV,\Icov\cV)$.
By Theorem \ref{shrink}, $\Icov\cV$ is shrinkable.
By Theorem \ref{main}, it remains to prove that $\Icov\cV$ is weakly Ramsey.

Let $\cV_1,\cV_2,\dots\nin\Icov\cV$ be pairwise disjoint with $\Union_{n\ge m}\cV_n\in\Icov\cV$ for all $m$.
For each $n$, let
$$\cU_n = \left\{\bigcap_{m\in I}V_m : I\sbst\N, |I|=n, (\forall m\in I)\ V_m\in\cV_m\right\}.$$

\bclm
$\cU_n\in\Icov\cV$.
\eclm
\bpf
As $\cV$ is closed under finite intersections, $\cU_n\sbst\cV$.
Assume that there is $B\in\cI$ not contained in any member of $\cU_n$.
Let $I=\{m : (\exists U\in\cV_m)\ B\sbst U\}$. Then $|I|<n$.
For each $m\in I$ choose $B_m\in\cI$ witnessing that
$\cV_m\nin\OI$. Then $B\cup\Union_{m\in I}B_m$ is not covered by any $U\in\Union_n\cV_n$, a contradiction.
\epf

Apply $\sone(\Icov\cV,\Icov\cV)$ to the sequence $\cU_n$, $n\in\N$,
to obtain elements $U_n\in\cU_n$ with $\sseq{U_n}\in\Icov\cV$.
Let $m_1$ be such that $V_{m_1}:=U_1\in\cV_{m_1}$. Inductively, for each $n>1$,
$U_n$ is an intersection of elements from $n$ many $\cV_m$-s, and thus
there are $m_n$ distinct from $m_1,\dots,m_{n-1}$, and an element $V_{m_n}\in\cV_{m_n}$,
such that $U_n\sbst V_{m_n}$. Then $\cA=\sseq{V_{m_n}}\in\Icov\cV$.
$\cA\sbst\Union_n\cV_n$, and $|\cA\cap \cV_n|\le 1$ for all $n$.
\epf

At the price of a slightly less combinatorial proof,
we can weaken the restriction of $\OI$-Lindel\"ofness substantially.

\bthm\label{unctble}
Assume that $X$ has a countable open $\cI$-cover.
Then the four items of Theorem \ref{sch++} are equivalent.
\ethm
\bpf
The proof is the same as that of Theorem \ref{RiswR}, but we argue directly in
some of its steps. We do this briefly.

$(1\Impl 4)$ By (1), $X$ is $\OI$-Lindel\"of, and the argument in the proof of
Theorem \ref{sch++} applies.

$(3\Impl 2)$ Let $\cU_1,\cU_2,\dots\nin\OI$ be disjoint with $\Union_n\cU_n\in\OI$. Set $\cU=\Union_n\cU_n$.
Define a coloring $c:[\cU]^2\to\{1,2\}$ by
$$c(U,V)=\begin{cases}
1 & (\exists n)\ U,V\in \cU_n\\
2 & \mbox{otherwise}
\end{cases}$$
By (3), there is a monochromatic $\cV\sbst\cU$ with $\cV\in\OI$.
It is easy to see that $\cV$ is as required in (2).

$(2\Impl 1)$ Use the premised $\sseq{U_n}\in\OI$:
Assume that $\cU_1,\cU_2,\dots\in\OI$. For each $n$, let
$$\cV_n=\{U_n\cap U : U\in\cU_n\}.$$
Now, $\Union_n\cV_n=\cU\in\Icov\cU$, and for each $n$, $\cV_n\nin\Icov\cU$.
By thinning out the sets $\cV_n$ if necessary, we may assume that they are disjoint.
By (2), there is $\cV\sbst\cU$ such that $\cV\in\OI$
and $|\cV\cap\cV_n|\le 1$ for all $n$.
\epf

For $T_1$ topological spaces, the assumption that $X$ has a countable open $\cI$-cover
can be simplified.

\blem
Let $\cI$ be an ideal on a $T_1$ space $X$. There is a countable $\cI$-cover of $X$ if,
and only if, there is a countable $D\sbst X$ such that $D\nin\cI$.
\elem
\bpf
$(\Rightarrow)$ Let $\cU$ be a countable $\cI$-cover of $X$. For each $U\in\cU$, pick $x_U\in X\sm U$.
Take $D=\{x_U : U\in\cU\}$.

$(\Leftarrow)$ $\cU=\{X\sm\{x\} : x\in D\}$ is a countable $\cI$-cover of $X$.
\epf

In particular, Scheepers's Theorem \ref{sch} is true \emph{for all $T_1$ spaces}:
It is trivially true for finite spaces, and in the remaining case there is a countably infinite subset.

In the case of $k$-covers, it suffices to assume that $X$ has a countable subset with noncompact closure.

\section{Back to van der Waerden's Theorem}

We reconsider van der Waerden's superfilter $\AP$ of all sets containing
arbitrarily long arithmetic progressions.

\bexm\label{FW}
Furstenberg and Weiss (unpublished) proved that $\AP\nrightarrow(\AP)^2_2$.
Using Theorem \ref{RiswR}, we can reproduce their observation by showing
that $\AP$ is not even \emph{weakly Ramsey}:
Let $A_1=\{1\}$, and for each $n>1$, let $m_n=2\max A_{n-1}$,
and $A_n=\{m_n+1,m_n+2,...,m_n+n\}$. For each $n$, $A_n\nin\AP$, and $\Union_nA_n\in\AP$.
But there is no arithmetic progression of length $3$ with at most one element in each $A_n$.
\eexm

Example \ref{FW} motivates us to look for a property which is weaker than being Ramsey
but still implies Ramsey's Theorem, and which is satisfied by $\AP$.
A natural candidate is available in the literature.

\bdfn[Baumgartner-Taylor \cite{BT78}]
$\cS\to\lceil\cS\rceil^n_k$ is the statement: For each $A\in\cS$ and each coloring $c:[A]^n\to\{1,2,\dots,k\}$,
there is $M\sbst A$ such that $M\in\cS$, and a partition of $M$ into finite pieces, such that
$c$ is constant on elements of $[M]^n$ containing at most one element from each piece.
\edfn

Any provable assertion of the form $\cS\to\lceil\cS\rceil^n_k$ with $\emptyset\neq\cS\sbst\roth$ and $n,k\ge 2$ is
an improvement of Ramsey's Theorem: Given a coloring of
$\N$, take $M\in\cS$ and a partition of $M$ into finite sets as promised by $\cS\to\lceil\cS\rceil^n_k$.
Then any choice of one element from each piece gives an infinite monochromatic set.
$\cS\to\lceil\cS\rceil^n_k$ also implies that $\cS$ is a superfilter.

\blem
For each upwards-closed $\emptyset\neq\cS\sbst\roth$:
\be
\itm If $\cS\to\lceil\cS\rceil^n_k$, $l\le n$, and $m\le k$, then  $\cS\to\lceil\cS\rceil^l_m$.
\itm For each $k$, $\cS\to\lceil\cS\rceil^1_k$ is equivalent to $\cS\to(\cS)^1_k$.
\ee
\elem
\bpf
(1) Given $c:[A]^l\to\{1,\dots,m\}$, define $f:[A]^n\to\{1,\dots,k\}$ by
letting $f(F)$ be the $c$-color of the $l$ smallest elements of $F$.
Use $\cS\to\lceil\cS\rceil^n_k$ to obtain $M\sbst A$ such that $M\in\cS$,
and a partition of $M$ into finite sets, such that sets with elements coming
from distinct pieces of $M$ all have the same $f$-color $i$.

For each $F\in[A]^l$ with elements coming from distinct pieces of $M$,
take arbitrary $n-l$ elements from other pieces of $M$, which are greater
than all elements of $F$ (this can be done since $M$ is infinite, and
the pieces are finite). Add these elements to $F$, to obtain $F'$.
Then $c(F)=f(F')=i$.

(2) Immediate from the definition.
\epf

\bdfn
A superfilter $\cS$ is
a \emph{$P$-point} if for all members $A_1\spst A_2\spst\dots$ of $\cS$, there is $A\in\cS$
such that $A\sm A_n$ is finite for all $n$.
\edfn

\bdfn[Scheepers \cite{coc1}]
$\sfin(\cS,\cS)$ is the statement: Whenever $S_1,S_2,\dots\in\cS$,
there are finite $F_n\sbst S_n$, $n\in\N$, such that $\Union_nF_n\in\cS$.
\edfn

\bthm\label{Ppt}
The following are equivalent for superfilters $\cS$:
\be
\itm $\cS$ is a $P$-point.
\itm $\sfin(\cS,\cS)$.
\itm For all disjoint $A_1,A_2,\dots$ with $\Union_{n\ge m}A_n\in\cS$
for all $m$, there is $A\sbst\Union_nA_n$ such that $A\in\cS$ and $A\cap A_n$ is finite for all $n$.
\itm For each partition $\N=\Union_nA_n$ with $\Union_{n\ge m}A_n\in\cS$
for all $m$, there is $A\in\cS$ such that $A\cap A_n$ is finite for all $n$.
\itm $\cS\to\lceil\cS\rceil^2_2$ and $\cS$ is shrinkable.
\itm $\cS\to\lceil\cS\rceil^n_k$ for all $n,k$, and $\cS$ is shrinkable.
\ee
\ethm
\bpf
$(1\Impl 2)$ Assume that $S_1,S_2,\dots\in\cS$.
For each $n$, let $A_n = \Union_{m\ge n} S_m$. By $(1)$,
there is $A\in\cS$ such that $A\sm A_n$ is finite for all $n$.
For each $n$, let $F_n = (A\cap S_n)\sm A_{n+1}$.
Let $B=A\cap\bigcap_n A_n$.
For each $n$, add at most finitely many elements of $B$
to $F_n$, in a way that $F_n$ remains finite, $F_n\sbst S_n$, and $\Union_n F_n\spst B$.
Then $A\sm\Union_n F_n$ is finite, and thus $\Union_n F_n\in\cS$.

$(2\Impl 3)$ apply $\sfin(\cS,\cS)$ to the sequence $\Union_{n\ge m}A_n$, $m\in\N$.

$(3\Impl 4)$ Trivial.

$(4\Impl 1)$ Assume that $B_1\spst B_2\spst\dots$ are members of $\cS$.
We may assume that $B_1=\N$. Let $A_0 =\bigcap_n B_n$. If $A_0\in\cS$ we are done, so assume that $A_0\nin\cS$.

For each $n$, let $A_n=B_n\sm B_{n+1}$.
$\N=A_0\cup\Union_nA_n$ is a partition of $\N$ as required in (3):
$\Union_nA_n\in\cS$ as $A_0\nin\cS$.
For each $n$, $\Union_{m\ge n}A_m=B_n\sm A_0\in\cS$, since $B_n\in\cS$.
Take $A\in\cS$ such that $A\cap A_n$ is finite for all $n$.
Then $A\sm B_n$ is finite for all $n$.

$(5\Impl 3)$ Consider disjoint $A_1,A_2,\dots$ with $\Union_{n\ge m}A_n\in\cS$
for all $m$. As $\cS$ is shrinkable, we may assume that $A_n\nin\cS$ for all $n$.
Let $D=\Union_{n}A_n$, and define a coloring $c:[D]^2\to\{1,2\}$ by
$$c(m,k)=\begin{cases}
1 & (\exists n)\ m,k\in A_n\\
2 & \mbox{otherwise}
\end{cases}$$
By $\cS\to\lceil\cS\rceil^2_2$, there is a partition $M=\Union_n F_n\sbst D$ into finite sets,
such that $M\in\cS$ and $c$ is constant on pairs of elements coming from different $F_n$-s.

Assume that these pairs have color $1$.
Fix $k\in F_1$, and $n$ such that $k\in A_n$.
For each $m\neq 1$ and each $i \in F_m$, $c(k,i)=1$ and thus $i\in A_n$, too.
But then each $l\in F_1$ has $c(i,l)=1$, and thus $l\in A_n$, too.
Thus, $M\sbst A_n$. As $M\in\cS$, we have that $A_n\in\cS$; a contradiction.
Thus, all pairs coming from different $F_n$-s, must come from different $A_n$-s.
Take $A=\Union_nF_n$.

$(1,3\Impl 6)$
Clearly, (3) implies that $\cS$ is shrinkable. We prove that $\cS\to\lceil\cS\rceil^d_k$ for all $d,k$,
by induction on $d$.

Let $\cS$ be a $P$-point superfilter,
$A\in\cS$, and $c:[A]^d\to\{1,\dots,k\}$.
The case $d=1$ follows from $\cS$ being a superfilter.

Induction step: Enumerate $A=\sseq{a_n}$.
Choose $A_{a_1}\sbst A\sm\{a_1\}$ and a color $i_{a_1}$ such that
$A_{a_1}\in\cS$, and a partition of $A_{a_1}$ into finite sets,
such that for each $F\in[A_{a_1}]^{d-1}$ with at most one element in each piece,
$c(\{a_1\}\cup F)=i_{a_1}$.
In a similar manner, choose inductively for each $n>1$
$A_{a_n}\sbst A_{a_{n-1}}\sm\{a_n\}$ and a color $i_{a_n}$ such that
$A_{a_n}\in\cS$, and a partition of $A_{a_n}$ into finite sets,
such that for each $F\in[A_{a_n}]^{d-1}$ with at most one element in each piece,
$c(\{a_n\}\cup F)=i_{a_n}$.

As $\cS$ is a $P$-point, there is $B\in\cS$ such that $B\sm A_{a_n}$
is finite for all $n$. Fix a color $i$ such that $C=\{n\in B : i_n=i\}\in\cS$.

Let $c_1=\min C$.
Inductively, for each $n>1$ choose $c_n\in C$ such that:
\be
\itm $c_n>c_{n-1}$;
\itm For each piece from the partitions of $A_{a_1},\dots,A_{a_n}$ which
intersects $[1,c_{n-1})$, $c_n$ is greater than all elements of that piece;
and
\itm $C\sm[1,c_n)\sbst A_{c_{n-1}}$.
\ee
As
$$C=\left(C\cap \Union_{n\in\N} [c_{2n},c_{2n+1})\right)\cup\left(C\cap \Union_{n\in\N}[c_{2n-1},c_{2n})\right),$$
there is $l\in\{0,1\}$ such that $M=C\cap \Union_n[c_{2n-l},c_{2n+1-l})\in\cS$.

Let $m_1<m_2<\dots<m_d$ be members of $M$ coming from distinct
intervals $[c_{2n-l},c_{2n+1-l})$. Let $n$ be minimal with $m_1<c_n$. Then
$$m_2,\dots,m_d\in C\sm[1,c_{n+1})\sbst A_{c_n}\sbst A_{m_1},$$
and $m_2,\dots,m_d$ come from distinct pieces of the partition of $A_{m_1}$.
Thus, $c(\{m_1,\dots,m_d\})=c(\{m_1\}\cup\{m_2,\dots,m_d\}) = i_{m_1}=i$.

$(6\Impl 5)$ Trivial.
\epf

The equivalence of (1) and (3) in the following corollary can be shown, using
a well known argument, to be the same as the equivalence of (i) and (iii) in Theorem 2.3 of
Baumgartner and Taylor \cite{BT78}.

\bcor
For ultrafilters $\cU$, the following are equivalent:
\be
\itm $\cU$ is a $P$-point.
\itm $\sfin(\cU,\cU)$.
\itm $\cU\to\lceil\cU\rceil^2_2$.
\itm $\cU\to\lceil\cU\rceil^n_k$ for all $n,k$.
\ee
\ecor
\bpf
Recall that ultrafilters are shrinkable.
\epf

\bdfn
A family $\cF$ of subsets of $\N$ \emph{generates} an upwards-closed family $\cS$ if
$\cF\sbst\cS$ and each element of $\cS$ contains an element of $\cF$.
An upwards-closed family $\cS\sbst\roth$ is \emph{compactly generated}
if there are upwards-closed families $\cF_1,\cF_2,\dots\sbst P(\N)$, each generated by finite subsets of $\N$,
such that $\cS=\bigcap_n\cF_n$.
\edfn

\bexm
$\roth$ is compactly generated: Take $\cF_n = [\N]^{\ge n}$, $n\in\N$.

$\AP$ is compactly generated: Let $\cF_n$ be the
family of all sets containing arithmetic progressions of length $n$.

Similarly, the \emph{Folkman-Rado-Sanders} superfilter \cite{Sanders}
of sets containing arbitrarily large finite subsets
together with all of their subset sums is compactly generated.
\eexm

Schur's Theorem \cite{Schur} states that if the natural numbers are colored
in finitely many colors, then there is a monochromatic
solution to the equation $x+y=z$.
Rado's Theorem \cite{Rado} extends Schur's Theorem to arbitrary
regular homogeneous systems of equations.
A homogeneous system of equations $Ax=0$ with integer coefficients
is \emph{regular} if the columns of $A$ can be partitioned into
sets $P_1,\dots,P_k$ such that $\sum_{v\in P_1}v=0$, and for each
$i>1$, each element of $P_i$ is a linear combination of elements
of $P_1\cup\dots\cup P_{i-1}$.

The family of all sets containing a solution to a regular homogeneous
system of equations is not a superfilter.
This problem can be solved by using the following operation on upwards-closed
families (see Proposition \ref{par} below).

\bdfn
For an upwards-closed family $\cF$ of subsets of $\N$ and $k\in\N$, $\Par_k(\cF)$
is the family of all $A\sbst\N$ such that for each partition
of $A$ into $k$ pieces, one of the pieces belongs to $\cF$.
$\Par(\cF)=\bigcap_k\Par_k(\cF)$.
\edfn

For upwards-closed families $\cF$, $\Par(\cF)\sbst\cF$, and
$\cF$  is a superfilter if, and only if, $\Par(\cF)=\cF$.


\blem\label{uc}
Assume that $\cF\sbst P(\N)$ is upwards-closed and generated by finite subsets of $\N$.
Then the same is true for $\Par_k(\cF)$, for all $k$.
\elem
\bpf
This is a reformulation of the compactness theorem for partitions, see Theorem 2.5 in \cite{ProtasovRT}.
\epf

Note that $\N\in\Par(\cF)$ if, and only if, $\Par(\cF)$ is nonempty.

\bprp\label{par}
Let $\cF$ be an upwards-closed family of subsets of $\N$.
Assume that $\cF$ does not contain any singleton, and $\N\in\Par(\cF)$. Then:
\be
\itm $\Par(\cF)$ is the maximal superfilter contained in $\cF$.
\itm If $\cF$ is compactly generated, then so is $\Par(\cF)$.
\ee
\eprp
\bpf
(1) It is easy to see that $\Par(\cF)$ is closed upwards.

Assume that $A\cup B\in\Par(\cF)$, and $A\cap B=\emptyset$.
If $A,B\nin\Par(\cF)$, then there are a partition of $A$ into $n$ pieces and a partition of
$B$ into $m$ pieces, such that none of the pieces belong to $\cF$.
But this yields a partition of $A\cup B$ into $n+m$ pieces, none of which
from $\cF$, that is, $A\cup B\nin\Par_{n+m}(\cF)$. A contradiction.

As $\Par(\cF)\sbst\cF$, there are no singletons in $\Par(\cF)$, and consequently
no finite sets.

If $\cS$ is any superfilter contained in $\cF$, then
$\cS=\Par(\cS)\sbst\Par(\cF)$.

(2) Assume that $\cF=\bigcap_n\cF_n$, with each $\cF_n$ upwards-closed and generated
by finite subsets. Replacing each $\cF_n$ by $\bigcap_{m\le n}\cF_m$, we may assume
that $\cF_1\spst\cF_2\spst\dots$.
It follows that for each $k$, $\Par_k(\bigcap_n\cF_n)=\bigcap_n\Par_k(\cF_n)$,
and thus
$$\Par(\cF)=\bigcap_{k\in\N}\Par_k(\cF)=\bigcap_{k\in\N}\Par_k(\bigcap_{n\in\N}\cF_n)=\bigcap_{k,n\in\N}\Par_k(\cF_n).$$
By Lemma \ref{uc}, each $\Par_k(\cF_n)$ is upwards-closed and generated by finite sets.
\epf

\bexm\label{schurex}
Let $\cF$ be the family of all subsets of $\N$
containing a solution to the equation $x+y=z$.
Let $\Par(x+y=z)=\Par(\cF)$.
Schur's Theorem tells that $\N\in\Par(x+y=z)$.
By Proposition \ref{par}, $\Par(x+y=z)$ is a compactly-generated
superfilter.
We can define similarly $\Par(Ax=0)$ for an arbitrary regular system
of homogeneous equations, and by Rado's Theorem have that
$\Par(Ax=0)$ is a compactly-generated
superfilter.
\eexm

We now state the main application of Theorem \ref{Ppt}.

\theoremstyle{theorem}
\newtheorem*{pithm}{Theorem $\pi$}

\begin{pithm}\label{BH++}
Assume that $\cS$ is a compactly-generated superfilter.
Then $\cS\to\lceil\cS\rceil^n_k$ for all $n,k$.
\end{pithm}
\bpf
By Theorem \ref{Ppt}, it suffices to show that $\sfin(\cS,\cS)$ holds.
Let $\cF_1,\cF_2,\dots\sbst P(\N)$ be upwards-closed and generated by
finite sets, such that $\cS=\bigcap_n\cF_n$.
Assume that $A_1,A_2,\dots\in\cS$. For each $n$, pick a finite $F_n\in\cF_n$
such that $F_n\sbst A_n$. Then $\Union_n F_n\in\cS$.
\epf

\stepcounter{thm}

\newcommand{\pit}{$\pi$} 
Theorem \pit{}
is a simultaneous improvement of the theorems of Ramsey, van der Waerden,
Schur, Rado, Folkman-Rado-Sanders, and many more.
In particular, we have the following.

\bcor\label{vDR}
$\AP\to\lceil\AP\rceil^n_k$, for all $n,k$. \qed
\ecor

Theorem \pit{} can be restated
as follows.

\bcor\label{BH++2}
Assume that $\cS$ is a superfilter compactly generated by $\cF_1,\cF_2,\dots$.
Then for all $r,k$, $A\in\cS$, $c:[A]^r\to \{1,\dots,k\}$, and $m_1<m_2<\dots$,
there are disjoint $F_n\in\cF_{m_n}$, $n\in\N$, such that
$\Union_nF_n\in\cS$, and $c$ is constant on sets with at most one element from each $F_n$.
\ecor
\bpf
We may assume that $\cF_1\spst\cF_2\spst\dots$.
Assume that $A\in\cS$ and $c:[A]^n\to\{1,2,\dots,k\}$.
Using Theorem \pit{}, take $M\sbst A$ such that $M\in\cS$, and a partition of $M$ into finite pieces, such that
$c$ is constant on sets containing at most one element from each piece.
$M$ contains some finite element of $\cF_{m_1}$.
Let $F_1$ be the union of as many pieces of $M$ as required so that $F_1$ contains
this element of $\cF_{m_1}$.
$M\sm F_1\in\cS$, and is partitioned by the remaining pieces, thus
we can take a union of finitely many of the remaining pieces,
$F_2$, containing some element of $\cF_{m_2}$, etc.

$\Union_nF_n$ contains an element of each $\cF_n$, and thus belongs to $\cS$.
\epf

\bexm
Consider Corollary \ref{BH++2} with $\cS=\AP$. Fix an arbitrarily quickly increasing
sequence $m_n$, and assume that we color an arbitrarily sparse $A\in\AP$.
Then each $F_n$ contains, and thus may be assumed to be, an arithmetic
progression of length $m_n$. The special case $A=\N$ is the main corollary in Bergelson and Hindman's 1989 paper \cite{BerHin89}.
\eexm

Bergelson and Hindman's proof in \cite{BerHin89}
shows that it suffices to assume that the colored set $A$ is an element
of a combinatorially large ultrafilter (see \cite{BerHin89}).
Elements of $\AP$ need not lie in a combinatorially large ultrafilter, and
we do not know a simple way to deduce Corollary \ref{vDR} (or \ref{BH++2})
from Bergelson and Hindman's Corollary, and not even from their much stronger Theorem 2.5 of \cite{BerHin89}.

\section{An additional application to topological selection principles}

Using Theorem \ref{Ppt} and arguments similar to those in the proof of Theorem \ref{sch++},
we also obtain the following Theorem \ref{OIPpt}. In the case that $\cI$ is the ideal of
finite sets ($\OI=\Omega$), the equivalence of (2) and (4) was proved by Just, Miller, Scheepers, and
Szeptycki in \cite{coc2}.
In the case that $\cI$ is the ideal of subsets of compact sets,
the equivalence of (2) and (4)
was proved by Di Maio, Ko\v{c}inac, and Meccariello in \cite{AppKI}.

\bthm\label{OIPpt}
Let $\cI$ be an ideal on $X$. If $X$ is $\OI$-Lindel\"of,
then the following are equivalent:
\be
\itm For all $\cU_1\spst\cU_2\spst\dots$ from $\OI$, there is $\cU\in\OI$ such that
$\cU\sm\cU_n$ is finite for all $n$.
\itm $\sfin(\OI,\OI)$.
\itm For all disjoint $\cU_1,\cU_2,\dots$ with $\Union_{n\ge m}\cU_n\in\OI$
for all $m$, there is $\cU\sbst\Union_n\cU_n$ such that $\cU\in\OI$ and $\cU\cap \cU_n$ is finite for all $n$.
\itm $\OI\to\lceil\OI\rceil^2_2$.
\itm $\OI\to\lceil\OI\rceil^n_k$ for all $n,k$.\qed
\ee
\ethm

Here too, by using direct arguments as in the proof of Theorem \ref{unctble},
``$X$ is $\OI$-Lindel\"of'' can be weakened to ``$X$ has a countable open $\cI$-cover'',
or equivalently for $T_1$ spaces, to ``there is a countable $D\sbst X$ such that $D\nin\cI$''.

\section{Final comments}

Mathias defines in \cite{Mathias77} \emph{happy families},
certain types of superfilters which were later named
\emph{selective} by Farah \cite{Farah}.
Farah points out in \cite{Farah} that every selective superfilter is Ramsey.
It is immediate that every selective superfilter is strongly Ramsey,
and arguments similar to those in the proof
of Lemma \ref{sRisR} show that every strongly Ramsey superfilter is selective.
Given Farah's observation, one can obtain a simpler proof of Lemma
\ref{sRisR}.

\forget
For an ideal $\cI$ on $\N$ (see before Definition \ref{afterideal}),
$\cI^+=\{A\sbst\N : A\nin\cI\}$ is called a \emph{coideal} (of the ideal $\cI$).
Superfilters are coideals, and vice versa.
Using this interpretation of superfilters,
one obtains a dual view of the results presented here.
We mention some of them.

Example \ref{RnotsR} is the coideal of
the ideal $\emptyset\x \mathrm{Fin}$ on $\N\x\N$,
consisting of all sets which are finite on each $x$-coordinate.
Also, the superfilter of somewhere dense sets from Example
\ref{wRnotR} is the coideal of the ideal of nowhere dense sets.

A Boolean algebra $B$ is \emph{$\sigma$-closed} if each sequence
$a_1\ge a_2\ge\dots$ in $B$ has a (nonzero) lower bound.
Thuemmel pointed out to us that, as can be easily verified,
a superfilter $\cI^+$ is shrinkable if,
and only if, the Boolean algebra $P(\N)/\cI$ is $\sigma$-closed.
In other words, a superfilter $\cS$ is shrinkable if, and only if,
for each sequence $A_1\spst A_2\spst\dots$ of element of $\cS$,
there is $A\in\cS$ such that for each $n$, $A\sm A_n\nin\cS$.

Thuemmel also observed that the weaker condition that $P(\N)/\cI$ is
$\alephes$-distributive suffices to make weakly Ramsey equivalent to
Ramsey for $\cI^+$.
\forgotten

Rec\l{}aw has informed us of his independent work with
Filip\'ow, Mro\.zek, and Szuca \cite{FMRS},
which contains related results, mainly of a descriptive set theoretic flavor.

In the topological results, considering from the start only countable covers
removes any restriction from the considered topological spaces.
For example, our results immediately apply to the corresponding families of countable \emph{Borel} covers,
since the Borel sets form a base for a topology on $X$.
A general study of countable Borel covers in the context of selection principles is available in \cite{CBC}.

Theorem \pit{} and its Corollary \ref{BH++2} should be viewed as a
simple way to lift one-dimensional Ramsey theoretic results to higher dimensions.
It does not generalize the Bergelson-Hindman Theorem from \cite{BerHin89},
but it extends it to cover additional classes of superfilters, and assumes less on the colored set.

\subsection*{Acknowledgments}
We thank Adi Jarden,  Jan Stary, Egbert Thuemmel, and Lyubomyr Zdomskyy for their useful comments
and suggestions, Andreas Blass for pointing out Berge's reference to us, and Vitaly Bergelson and Neil Hindman
for their encouraging reaction to these results, when presented in the conference \emph{Ultramath 2008}
(Pisa, Italy, June 2008). There is no doubt that these results owe much to the pioneering works
of Scheepers and his followers, some of which are mentioned in the bibliography, and in the references therein.

\ed